\documentstyle[12pt]{article}
\newtheorem{proposition}{Proposition}
\newtheorem{corollary}{Corollary}

\newtheorem{remark}{Remark}
\newcommand{\ot}{\otimes}

\newcommand\De{\Delta}
\newcommand\ar{{\bf A}(S)\,}
\newcommand\vv{{\bf V}^{\ot 2}}

\newcommand\RR{\cal R}

\newcommand{\wpp}{\wedge_+}
\newcommand{\wm}{\wedge_-}
\newcommand{\Ppm}{P_{\pm}}
\newcommand{\wpm}{\wedge_{\pm}}
\newcommand{\uq}{U_q(\gggg)}
\newcommand{\gggg}{g}

\newcommand{\Pm}{P_-(t)}
\newcommand{\rr}{\rho^{\ot 2}}
\def\adots{\mathinner{\mkern2mu\raise1pt\hbox{.}
\mkern3mu\raise4pt\hbox{.}\mkern1mu\raise7pt\hbox{.}}}
\newcommand{\Fun}{{\rm Fun}}
\newcommand{\llangle}{\langle\langle}
\newcommand{\rrangle}{\rangle\rangle}

\title{Dual quasitriangular structures related to the Temperley--Lieb
algebra}

\author{P.~Akueson and D.~Gurevich,\\
ISTV, Universit\'e de Valenciennes
59304 Valenciennes, France \\
e-mail:  akueson@univ-valenciennes.fr, gurevich@univ-valenciennes.fr}

\begin{document}
\maketitle 
\begin{abstract}
We consider nonquasiclassical solutions to the quantum Yang--Baxter equation
and the corresponding quantum cogroups $\Fun(SL(S))$ constructed earlier in
\cite{G}. We give a criterion of the existence of a dual quasitriangular
structure in the algebra $\Fun(SL(S))$ and describe a large class of such
objects related to the Temperley--Lieb algebra satisfying this criterion. We
show also that this dual quasitriangular structure is in some sense
nondegenerate.

 \end{abstract}

\section{Introduction}

It became clear after the works of one of the authors \cite{G}
that besides the well-known deformational (or
quasiclassical) solutions to the quantum Yang--Baxter
equation (QYBE) there exists a lot of other solutions that
differ drastically from the former ones. Let us explain this
in more detail. Let ${\bf V}$ be a linear space over the
field $K={\bf C} $
 or $ R$. We call a {\em Yang--Baxter operator}
a solution $S:\vv\to\vv$ to the QYBE
$$S^{12}S^{23}S^{12}=S^{23}S^{12}S^{23},\;S^{12}=S\ot
{\rm id},\; S^{23}={\rm id}\ot S.$$ A Yang--Baxter
operator satisfying a second degree equation $$({\rm
id}+S)(q\,{\rm id}-S)=0$$ was called in \cite{G} a
{\em Hecke symmetry}. The quantum parameter $q\in K$
is assumed to be generic.

It is natural to associate to a Hecke symmetry two algebras defined as
follows
$$\wpp({\bf V})=T({\bf V})/\{{\rm Im}(q\,{\rm id}-S)\},\;\wm({\bf
V})=T({\bf V})/\{{\rm Im}({\rm id}+S)\}.$$
They are q-counterparts of the symmetric and skew-symmetric algebras of the
space ${\bf V}$, respectively.

Let us denote $\wpm^l({\bf V})$ the degree $l$ homogeneous component of the
algebra $\wpm({\bf V})$. It was shown in \cite{G} that the {\em Poincar\'e
series}
$$\Ppm(t)=\sum {\rm dim}\wpm^l({\bf V}) t^l$$
of the algebras $\wpm({\bf V})$ for a generic $q$ satisfy the standard
relation
$$P_+(t) P_-(-t)=1.$$
Moreover, if the series $\Pm$ is a polynomial with 
leading coefficient 1, it
is reciprocal. Hecke symmetries of such type and the corresponding linear
spaces ${\bf V}$ are called {\em even}.

A particular case of an even Hecke symmetry is provided by the quantum groups
$U_q(sl(n))$: the operator $S=\sigma\rr(\RR)$ (where ${\cal R}$ is the
corresponding universal R-matrix, $\sigma$ is the flip and $\rho:sl(n)\to{\rm
End}({\bf V})$ is the fundamental vector representation) is just such a type
of solution to the QYBE. In fact, we have a family of operators $S_q$ and we
recover the standard flip $\sigma$ for $q=1$. Namely, in this sense we call
such Yang--Baxter operators (and all related objects) {\em deformational} or
{\em quasiclassical}.

The Poincar\'e series for Hecke symmetries of such type coincide with the
classical ones. Thus, in this case we have $P_-(t)=(1+t)^n$ with $n=\dim\,
{\bf V}$ and, therefore, $P_+(t)=(1-t)^{-n}$.

However, this is no longer true in general case. As shown in \cite{G}, for
any $n=\dim\,{\bf V}$ and any integer $p,\;2\leq p\leq n$ there exists a
nonempty family of nontrivial even Hecke symmetries such that ${\rm deg}\,
P_-(t)=p$. The integer $p$ is called the {\em rank} of the even Hecke symmetry
$S$ or of the corresponding space ${\bf V}$.

The classification problem of all even Hecke symmetries of a given rank $p$ is
still open. However, all such symmetries of rank $p=2$ are completely
classified.

Let us observe that the case $p=2$ is related to the Temperley--Lieb (TL)
algebra, since the projectors $P^i_-:{\bf V}^{\ot m}\to {\bf V}^{\ot m},\; 1
\leq i\leq m-1$ defined by
$$P^i_-={(q\,{\rm id}-S^{i,i+1})\over (q+1)}$$
where $S^{i,i+1}$ is the operator $S$ acting onto the $i$-th and $(i+1)$-th
components of ${\bf V}^{\ot m}$, generate a TL algebra. Let us recall that a
TL algebra is the algebra generated by $\,t_i,\,\,1\leq i\leq n-1\,$ with the
following relations:
$$t_i^2=t_i,\,\,t_i\,t_{i\pm 1}\,t_i=\lambda t_i\,\,,
t_i\,t_j=t_j\,t_i\,\,\,\,\vert i-j\vert\,>1.$$
(In the case under consideration $\lambda=q(1+q)^{-2}$.) In what follows even
rank 2 Hecke symmetries will be called symmetries of {\em TL type}.

It is possible to assign to any YB operator $S$ the famous ``RTT=TTR"
algebra. In the sequel it is called the {\em quantum matrix algebra} and
denoted by $\ar$. If $S$ is an even Hecke symmetry of the TL type, in this
algebra there exists a so-called {\em quantum determinant} (it was introduced
in \cite{G}). If it is a central element, it is possible to define a {\em
quantum cogroup} $\Fun(SL(S))$ looking like the famous quantum function algebra
$\Fun_q(SL(n))$ (in \cite{G} this algebra was called a quantum
group)\footnote{Let us note that a subclass of objects of such type was
independently introduced in \cite{DL}.}. Let us note that these quantum
cogroups $\Fun(SL(S))$ possess Hopf algebra structures.

However, until now no corepresentation theory of such nonquasiclas\-sical
quantum cogroups has been constructed. In particular, it is not clear whether
any finite-dimensional $\Fun(SL(S))$-comodule is semisimple for a generic
$q$. Nevertheless this problem seems to be of great interest, since the
nonquasiclassical solutions to the QYBE provides us a new type symmetries,
which differ drastically from classical or supersymmetries. (The simplest
models possessing such a symmetry of new type corresponding to an involutive
$S$, namely, a ``nonquasiclassical harmonic oscillator" was considered in
\cite{GRZ}. Let us observe that the partition functions of such models can be
expressed in terms of the Poincar\'e series corresponding to the initial
symmetry.)

The present paper is the first in a series aimed at a better understanding 
the structure of such nonquasiclassical symmetries. More precisely, we
discuss here two problems: first, whether the quantum cogroups $\Fun(SL(S))$
have quasi\-triangular structure and second, what is an explicit description
of their dual objects?

It is well known that the notion of a quasitriangular structure was
introduced by V.~Drinfeld. In fact this notion was motivated by the quantum
groups $\uq$. These objects have an explicit description due to V.~Drinfeld
and M.~Jimbo in terms of deformed Cartan--Weyl system $\{H_{\alpha},\,X_{\pm
\alpha}\}$ (cf., i.e., \cite{CP}). Thus this construction allows one to
develop a representation theory of quantum groups.

In the nonquasiclassical case under consideration such a description does not
exist. And the problem of an appropriate description of objects dual to the
quantum cogroups $\Fun(SL(S))$ is of great interest. (The duality in the
present paper is understood in the algebraic sense, i.e., all dual objects
are {\em restricted}).

We attack this problem here by means of the so-called {\em canonical
pairing}. Such a pairing can be defined on any algebra $\ar$ (an algebra
$\ar$ equipped with such a pairing is called {\em dual quasitriangular}).
Nevertheless, only under some additional conditions this pairing can be
descended to the quantum cogroup $\Fun(SL(S))$. We show here that this
condition is satisfied for the quantum cogroup $\Fun(SL(S))$ related to TL
algebras.

Moreover, we show that in this case the canonical pairing is nondegenerate
when restricted to the span of the generators of this algebra. This is the
main difference between the quasiclassical and nonquasiclassical cases: in
the former case this pairing is degenerate. This is a reason why we cannot
introduce an object dual to $\Fun_q(SL(n))$ by means of this pairing.
Finally, following the paper \cite{RTF}, in this case we must introduce an
additional pairing (associated in a similar way to the YB operator $S^{-1}$).
And the above mentioned deformed Cartan--Weyl basis in $\uq$ can be
constructed by means of both pairings (for this construction the reader is
referred to \cite{RTF}, cf. also Remark 1).

As for the nonquasiclassical case we can equip the basic space ${\bf V}$
(using nondegeneracy of the canonical pairing in the mentioned sense) with a
structure of a $\Fun(SL(S))$-module. Moreover, we can equip any tensor power
of the space ${\bf V}$ with such a structure. Finally, we get a new tool to
study tensor categories generated by such spaces.

The paper is organized as follows. In the next Section we introduce the
notion of a dual quasitriangular structure and describe the one connected
with the quantum matrix algebras in terms of the canonical pairing. In
Section 3 we give the condition ensuring the existence of such a structure on
the algebra $\Fun(SL(S))$ mentioned above and in Section 4 we show that this
condition is satisfied for a large family of such algebras related to the
TL algebras. We conclude the paper with a proof of the nondegeneracy of the
canonical pairing (in the above sense) for algebras from this family (Section
5) and with a discussion of a hypothetical representation theory for the
algebra $\Fun(SL(S))$ (Section 6).

\section{Dual quasitriangular structure}
The notion of a dual quasitriangular bialgebra (in particular, a Hopf
algebra) was introduced by Sh.~Majid (see \cite{M} and the references
therein) as a dualization of the notion of a quasitriangular bialgebra due to
V.~Drinfeld. By definition, a dual quasitriangular bialgebra is a bialgebra
equipped with some pairing similar to the one  defined on the cogroups
$\Fun_q(G)$ by means of the quantum universal R-matrix ${\cal R}$:
$$a\ot b\to\llangle a,b\rrangle=\langle a\ot b,\,{\cal R}\rangle
\,\, a,\, b\in\Fun_q(G).$$
Here the pairing $\langle\,\,,\,\,\rangle $ is that between $\Fun_q(G)$ and
$\uq$ extended to their tensor powers.

More precisely, one says that a bialgebra (or a Hopf algebra) ${\cal A}$ is
equipped with a {\em dual quasitriangular structure} and it is called a {\em
dual quasitriangular bialgebra} (or {\em Hopf algebra}), if it is endowed
with a pairing
$$
\llangle\,\,,\,\,\rrangle :{\cal A}^{\ot 2}\to K$$
satisfying the following axioms
\begin{eqnarray*}
{\rm (i)}&&\quad\llangle a,\,bc\rrangle=\llangle a_{(1)},\,c\rangle
\rangle\llangle a_{(2)},\,b\rrangle,\\
{\rm (ii)}&&\quad\llangle ab,\, c\rrangle=\llangle a,\,{c_{(1)}}
\rrangle\llangle b,\,{c_{(2)}}\rrangle,\\
{\rm (iii)}&&\quad\llangle a_{(1)},\,b_{(1)}\rrangle a_{(2)}b_{(2)}=b_{(1)}
a_{(1)}\llangle a_{(2)},\,{b_{(2)}}\rrangle,\\
{\rm (iv)}&&\quad\llangle 1,\, a\rrangle=\varepsilon (a)=\langle
\langle a,\,1\rrangle
\end{eqnarray*}
for all $\,a\,,b\,,c\,\,\in{\cal A}$, where $\varepsilon:{\cal A}\to K$ is
the counit of ${\cal A}\,$ and $\De: {\cal A}\to {\cal A}^{\ot 2},\,\De (a)=
a_{(1)}\ot a_{(2)}$ is the coproduct. If $\cal A$ is a Hopf algebra and
$\gamma:\cal A\to\cal A$ is its antipode, we impose a complementary axiom
$$
{\rm (v)}\quad\llangle a,\,b\rrangle=\llangle\gamma(a),\,\gamma(b)\rrangle.
$$
(If the pairing is invertible in sense of \cite{M}, p.~48 the axioms
(i)--(iii) imply those (iv) and (v), cf \cite{M}.)

Let us note that the axioms (i), (ii), (iv), (v) mean that the product
(resp., the coproduct, the unit, the counit, the antipode) of the algebra
${\cal A}$ is dual to the coproduct (resp. the product, the counit, the unit,
the antipode) of the algebra ${\cal A}^{op}$ where ${\cal A}^{op}$ denotes as
usually the bialgebra ${\cal A}$ whose product is replaced by the opposite
one. So, in fact, we have the pairing of bialgebras (Hopf algebras)
\begin{equation}
\llangle\,\,,\,\,\rrangle :{\cal A}\ot {\cal A}^{op}\to K.
\end{equation}

In some sense the notion of a dual quantum bialgebra is more fundamental than
that of quasitriangular one for the following reason. It is well known that
the most popular construction of a quasitriangular Hopf algebra is given by
the famous Drinfeld--Jimbo quantum group $\uq$. Usually it is introduced by
means of the Cartan--Weyl generators $\{H_{\alpha}, X_{\pm\alpha}\}$ and
certain relations between them which are quantum (or q-) analogues of the
ordinary ones. However, this approach is valid only in the quasiclassical
case.

In the general case, including  the nonquasiclassical objects, we should first
introduce dual quasitriangular bialgebras (or Hopf algebras) and only after
that we can proceed to introduce their dual objects. Moreover, an explicit
description of the latter objects depends on the properties of the canonical
pairing and they are not similar in quasiclassical and nonquasiclassical
cases.

Let us describe now a regular way to introduce the dual quasitriangular
bialgebras (and Hopf algebras) associated to the YB operators discussed
above. Let ${\bf V}$ be a linear space equipped with a nontrivial
Yang--Baxter operator $S:\vv\to\vv$. Let us fix a basis $\{x_i\}$ in ${\bf
V}$ and denote by $S_{ij}^{kl}$ the entries of the operator $S\,\,\,(S(x_i\ot
x_j)=S_{ij}^{kl} x_k\ot x_l)$. From here on summation on repeated indices is
assumed.

Let us consider a matrix $t$ with entries $t_k^l ,\,1\leq k,l\leq n={\rm dim}
\,{\bf V}$. The bialgebra $\ar$ of {\em quantum matrices} associated to $S\,$ is
defined as the algebra generated by $1\,$ and $n^2\,$ indeterminates $\{t_k^l\}
\,$, satisfying the following relations
$$S(t\ot t)=(t\ot t)S,\,\,{\rm or \,\,in\,\,a\,\, basis\,\,form}\,\,
S_{ij}^{mn}t_m^p\,t_n^q=t_i^u\,t_j^v S_{uv}^{pq}.$$
This algebra possesses a bialgebra structure, being equipped with the
co-matrix coproduct $\De(1)=1\ot1,\,\De(t_i^j)=t_i^p\ot t_p^j$ and the counit
$\varepsilon (1)=1,\,\varepsilon(t_i^j)=\delta_i^j$.

This is just the famous ``RTT=TTR" bialgebra introduced in \cite{RTF}. Let us fix
$c\in K,\,c\not=0$, and equip this algebra with a dual quasitriangular structure
by setting
$$\llangle 1,\,t_i^k\rrangle _c=\delta_i^k=\llangle t_i^k,\,1\rrangle _c\,\,
{\rm and}\,\,\llangle t_i^k,\,t_j^l\rrangle _c=c\,S_{ji}^{kl}$$
and extending the pairing to the whole $\ar^{\ot 2}$ by using the above
axioms (i), (ii) and (iv).

We leave to the reader to check that this extension is well defined (here it
is precisely the QYBE that plays the crucial role) and, moreover, the axiom
(iii) is satisfied as well.

The pairing $\llangle\,\,,\,\,\rrangle _c$ will be called {\em canonical}.

Let us remark that such a canonical pairing is usually considered with $c=1$.
However, we need this complementary ``degree of freedom" to make the pairing
$\llangle\,\,,\,\,\rrangle _c$ compatible with the equation $\det\,t=1$ (cf.
below). Another (but equivalent) way consists in replacing the Hecke symmetry
$S$ by $cS$. We drop the index $c$ if $c=1$.

Thus, the bialgebra $\ar$ can be {\em canonically} equipped with a dual
quasitriangular structure.

Nevertheless, this bialgebra does not possess any Hopf algebra structure
since no antipode is defined in it. To get a Hopf algebra, we must either
impose the complementary equation $\det\,t=1$, i.e., pass to the quotient
of the algebra $\ar$ by the ideal generated by the element $\det\,t-1$
(assuming the determinant $\det\,t$ to be well defined)\footnote{In the
sequel we will restrict ourselves to quantum cogroups of $SL$ type whose
construction was suggested in \cite{G}. Quantum cogroups of $SO$ or $Sp$ type
and the corresponding dual quasitriangular structures will be discussed
elsewhere.} or add to the algebra a new generator $\det^{-1}$. In the latter
case we obtain Hopf algebras (quantum cogroup) of $GL$ type
(cf. \cite{G}). 

In any case it is necessary to check that the above dual quasitriangular
structure on the algebra $\ar$ can be transferred to the final quantum
cogroup. As for standard quantum function algebras $\Fun_q(G)$ dual to the
quantum groups $\uq$ (for a classical simple Lie algebra $\gggg$) this
follows automatically by duality. However, in general this is no longer true.
In Section 3 we will give a necessary and sufficient condition ensuring the
existence of a dual quasitriangular structure on a $SL$ type quantum cogroup.

Thus, it is possible to associate a dual quasitriangular bialgebra to any YB
operator and a dual quasitriangular Hopf algebra (of $SL$ type) to some of
them. These algebras look like the function algebras $\Fun(G)$ on a ordinary
(semi) group $G$. This means that we can equip the space ${\bf V}$ with a
(right to be concrete) comodule structure over the algebra $\Fun(SL(S))$ by
$$\De: {\bf V}\to {\bf V}\ot\ar,\,\,\De(x_i)=x_j\ot t_i^j.$$
Therefore any tensor power of the space $\bf V$ also becomes a right
$\ar$-comodule.

However, these powers are not in general irreducible as comodules over the
above coalgebra. Unfortunately, up to now no corepresentation theory of
quantum cogroups in nonquasiclassical cases  has been constructed yet (such a
hypothetical theory in the case connected to the TL algebra is discussed in
Section 6). It is worth saying that even in the quasiclassical case it is
possible to use quantum cogroups $\Fun_q(G)$ instead of the quantum groups
$\uq$ themselves. However, technically it is more convenient to work with the
latter objects.

In the nonquasiclassical case an interesting problem arises: what is an
appropriate description of the objects dual to the bialgebras $\ar$ or of
their quotient of the $SL$ type. If the canonical pairing (1) is
nondegenerate, we can consider the bialgebra $\ar^{op}$ as the dual object to
that $\ar$ (and similarly for their quotients of the $SL$ type).

This is just (conjecturally) the case of the Hecke symmetries of TL type. We
show that (at least for a large family of such symmetries) the canonical
pairing, being restricted to the space ${\bf T}={\rm Span}(t_i^j)$, is
nondegenerate for a generic $q$. Nevertheless, this weak version of
nondegeneracy is sufficient to equip the initial space ${\bf V}$ with a
structure of a left $\ar^{op}$-module (and therefore, with that of a right
$\ar$-module).

Having in mind the usual procedure (cf., i.e., \cite{M}) we put
$$t_i^j\triangleright x_k=x_m\llangle t_k^m ,t_i^j\rrangle _c=cS_{ik}^{mj}
x_m,$$
where $t_i^j\triangleright x_k$ denotes the result of applying the element
$t_i^j\in\ar^{op}$ to $x_k\in{\bf V}$.

\begin{remark}{\em Let us observe that if the canonical pairing is
degenerate, then the above action is still well defined, but ${\bf V}$ becomes
reducible as an $\ar^{op}$-module since it contains the $\ar^{op}$-module
{\rm Im}($\,\triangleright$), where $\triangleright:\ar^{op}\ot {\bf V}\to
{\bf V}$ is the above map and this module is a proper submodule in ${\bf V}$.

It is just the case related to the quantum groups $\uq$. This is a
reasons why one needs a complementary pairing. More precisely, let us
introduce (following \cite{RTF}) two sets of generators $(L^+)_i^j$ and $(L^
-)_i^j$ and define the pairing between the spaces ${\bf L}^+={\rm Span}((L^+
)_i^j)$, ${\bf L}^-={\rm Span}((L^-)_i^j)$ and ${\bf T}={\rm Span}(t_i^j)$ as
follows
$$\llangle t_i^k,\,(L^+)_j^l\rrangle=S_{ji}^{kl},\,\,\llangle t_i^k
,\,(L^-)_j^l\rrangle=(S^{-1})_{ji}^{kl}.$$
In fact, in this way we have introduced a pairing between the spaces ${\bf L}
^+\oplus {\bf L}^-$ and ${\bf T}$. Of course, this pairing is degenerate on
${\bf L}^+\oplus {\bf L}^-$ but it becomes nondegenerate on ${\bf T}$. There
exists a natural way to extend the above pairing up to that $({\bf L}^+\oplus
{\bf L}^-)\ot\ar\to K$, cf. \cite{RTF}. Thus, the space ${\bf L}^+\oplus {\bf
L}^-$ is embedded into the algebra $\ar^*$ dual to $\ar$. The subalgebra of
$\ar^*$ generated by 1 and the space ${\bf L}^+\oplus {\bf L}^-$ is called in
\cite{RTF} the {\em algebra of regular functions} on $\ar$. (Moreover, in
\cite{RTF} the elements $(L^+)_i^j$ are expressed in terms of the generators
of the quantum groups $\uq$.) In a similar way we can define such an algebra
in the nonquasiclassical case under consideration, but since the canonical
pairing is nondegenerate for a generic $q$, we  restrict ourselves to the
generators $(L^+)_i^j$.}\end{remark}

\section{Dual quasitriangular algebras related to the even Hecke symmetries}
First, we recall some facts about the cogroups $\Fun(SL(S))$ introduced in
\cite{G}.

Let us fix an even Hecke symmetry $S:{\bf V}^{\ot 2}\to {\bf V}^{\ot 2}$ of
rank $p\geq 2$. Let us denote by $P_-^{(p)}$ the projector of ${\bf V}^{\ot
p}$ onto its skew-symmetric component $\wm^p({\bf V})$ (an explicit form of
this projector is given in \cite{G}). Then by definition ${\rm dim\, Im}\,P_-
^{(p)}=1$ and (assuming a base $\{x_i\}\in {\bf V}^{\ot p}$ to be fixed)
$$P_-^{(p)}x_{i_1}x_{i_2}\dots x_{i_p}=u_{i_1i_2\dots i_p}v^{j_1j_2\dots j_p}
x_{j_1}x_{j_2}\dots x_{j_p}$$
with $u_{i_1i_2\dots i_p}v^{i_1i_2\dots i_p}=1$ (hereafter we drop the sign
$\ot$).

The tensors $U=(u_{i_1i_2\dots i_p})$ and $V=(v^{j_1j_2\dots j_p})$ are
quantum analogue of the Levi--Civita ones.

Let us consider the bialgebra $\ar$ corresponding to the given even Hecke
symmetry $S$ and introduce a distinguished element in it
$$\det\,t=u_{i_1i_2\dots i_p}t_{j_1}^{i_1}\dots t_{j_p}^{i_p} v^{j_1j_2\dots
j_p}.$$
In \cite{G} it was shown that this element is group-like, i.e.,
$$\De\,(\det\,t)=\det\,t\ot\det\,t.$$
It was called a {\em quantum determinant}.

Under the additional condition that this determinant is central (in general
this is not so), we introduce an analogue $\Fun(SL(S))$ of the quantum
functional algebra $\Fun_q(SL(n))$ as the quotient algebra of $\ar$ over the
ideal generated by the element $\det\,t-1$. This quotient inherits a
bialgebra structure but, moreover, it possesses a Hopf structure (for an
explicit description of the antipode, the reader is referred to \cite{G}).
Our intermediate aim is to study whether it is possible to equip the algebra
$\Fun(SL(S))$ with a dual quasitriangular structure?

It is evident that the dual quasitriangular structure on $\ar$ defined above
can be descended to $\Fun(SL(S))$ iff
\begin{equation}
\llangle\det\,t,\,a\rrangle _c=\varepsilon(a)=\llangle a,\,\det\,t\rrangle _c
\end{equation}
for any $a\in\ar$. Using the fact that the quantum determinant is a
group-like element, it is possible to show that these relations are valid for
any $a$ if they are true for $a=t_i^j$. (As for $a=1$, relation (2) follows
immediately from $u_{i_1i_2\dots i_p}v^{i_1i_2\dots i_p}=1$.) Moreover, we
have
\begin{proposition} We have the following relations
\begin{equation}
\llangle t_k^l,\,\det\,t\rrangle _c=c^p(-1)^{p-1}qp_qM_k^l,\,\llangle\det\,t,
\,t_k^l\rrangle _c=c^p(-1)^{p-1}qp_qN_k^l,
\end{equation}
where $M_k^l=u_{i_1i_2\dots i_{p-1}k}v^{li_1i_2\dots i_{p-1}}$, $N_k^l=u_{ki
_1i_2\dots i_{p-1}}v^{i_1i_2\dots i_{p-1}l}$ and $p_q=1+q+\dots +q^{p-1}$.
{\em (}Let us note that the operators $M=(M_k^l)$ and $N=(N_k^l)$ have been
introduced in {\em \cite{G}, p.~816.)}
\end{proposition}

{\it Proof.} By axiom (i) we have
$$
\llangle t_k^l,\,\det\,t\rrangle _c=c^p u_{i_1i_2\dots i_p}v^{j_1j_2\dots
j_p}\llangle t_{k}^{m_1},t_{j_p}^{i_p}\rrangle\llangle t_{m_1}^{m_2},
t_{j_{p-1}}^{i_{p-1}}\rrangle\dots
$$
$$
\llangle t_{m_{p-1}}^{l},t_{j_1}^{i_1}\rrangle=c^pu_{i_1i_2\dots i_p}v^
{j_1j_2\dots j_p}S_{j_pk}^{m_1i_p}S_{j_{p-1}m_1}^{m_2i_{p-1}}\dots S_{j_1
m_{p-1}}^{li_1}.
$$
The term $u_{i_1i_2\dots i_p}v^{j_1j_2\dots j_p}S_{j_pk}^{m_1i_p}S_{j_{p-1}
m_1}^{m_2i_{p-1}}\dots S_{j_1m_{p-1}}^{li_1}$ was found in \cite{G} while
commuting the elements $V=v^{j_1j_2\dots j_p}x_{j_1}x_{j_2}\dots x_{j_p}$ and
$x_k$ (cf. Proposition 5.7 from \cite{G}) and is equal to $(-1)^{p-1}qp_qM_k
^l$. This proves the first equality. The second one can be proved in the same
way using the commutation law of the elements $x_k$ and $V=v^{j_1j_2\dots
j_p}x_{j_1}x_{j_2}\dots x_{j_p}$. \ $\Box$

\begin{corollary} Equations {\em (2)} can be satisfied for some $c\in K$ iff
the operators $M$ and $N$ are scalar {\em (}this property is equivalent by
virtue of Proposition {\em 5.9} from {\em \cite{G}} to the quantum
determinant being central{\em )} and, moreover, $M=N$. More precisely, if
$M=m\,{\rm id},\, N=n\,{\rm id},\,m,\,n\in K$, and $m=n$ we can satisfy the
relations
$$\llangle t_k^l,\,\det\,t\rrangle _c=\delta_k^l=\llangle\det\,t,\,t_k^l
\rrangle _c$$
by putting $c^p=(-1)^{p-1}q^{-1}p_q^{-1}m^{-1}$.
\end{corollary}

Let us note that the operators $M$ and $N$ satisfy the relation $MN=q^{p-1}p_
q^{-2}\, {\rm id}$ (cf. \cite{G}). Thus, if $M=m\, {\rm id}$, $N=n\,{\rm
id}$, the relation $M=N$ is equivalent to
\begin{equation}
m^2=q^{p-1}p_q^{-2}.
\end{equation}
Thus, we have reduced the problem of describing the quantum cogroups
$\Fun(SL(S))$ allowing a dual quasitriangular structure to the classification
problem of all even Hecke symmetries such that the corresponding operator
$M$ is scalar, $M=m\, {\rm id}$, with $m$ satisfying (4). In the next section
we will consider this problem for Hecke symmetries of TL type.
\begin{remark}
{\em Let us observe that if the operators $M$ and $N$ are not scalar, one
cannot define the algebra $\Fun(SL(S))$, but it is possible to define the
algebra $\Fun(GL(S))$ by introducing a new generator $\det^{-1}$ satisfying
the relations $\det^{-1}\,\det\,t=1$, and the commutation law of $\det^{-1}$
with other generators arising from this relation (cf. \cite{G}). Moreover, it
is possible to extend the canonical pairing up to that defined on
$\Fun(GL(S))$ by setting
$$\llangle t_i^p,\,{\det}^{-1}\rrangle\llangle t_p^j,\,\det\,t\rrangle=\delta
_i^j,\,\llangle{\det}^{-1},\,t_i^p\rrangle\llangle\det\,t,\,t_p^j\rrangle=
\delta_i^j.$$
The details are left to the reader.}
\end{remark}
\begin{remark}
{\em Let us observe that if an even Hecke symmetry is of TL type and $M$ is
scalar, then $M=N$ since in this case $M=VU$ and $N=UV$ (cf. Section 4).
Therefore, if the algebra $\Fun(SL(S))$ is well defined (i.e., the
corresponding quantum determinant is central) it automatically has a
canonical dual quasitriangular structure. It is not clear whether there
exists a Hecke symmetry of rank $p>2$ such that the algebra $\Fun(SL(S))$ is
well defined (i.e., $M=m\,{\rm id}$) but the factor $m$ does not satisfy the
relation (4) and therefore the corresponding canonical pairing is not
compatible with the equation $\det\,t=1$.}
\end{remark}

\section{The TL algebra case}
Now let us consider the case related to TL algebras. In this case it is
possible to give an exhausting classification of the corresponding Hecke
symmetries.

Indeed, it is easy to see (cf. \cite{G}) that any even Hecke symmetry of TL
type can be expressed by means of the Levi--Civita tensors $U=(u_{ij})$ and
$V=(v^{kl})$ in the following way
$$S_{ij}^{kl}=q\delta_i^k\delta_j^l-(1+q)u_{ij}v^{kl}.$$
Then the QYBE and the Hecke second degree relation are equivalent to the
system
\begin{equation}
{\rm tr}\,UV^t=1,\,\,\,UVU^{t}V^{t}=q(1+q)^{-2}\,{\rm id}.
\end{equation}
Hereafter $U\to U^t$ is the transposition operator. Thus, ${\rm tr}\,UV^t=
u_{ij}v^{ij}$.

Introducing the matrix $Z=(1+q)VU^t\,\,(z_i^j=(1+q)v^{jk}u_{ik})$ and using
the fact that the second relation of (5) can be represented in form $V^tU^tVU
=q(1+q)^{-2}\,{\rm id}$, we can reduce the relations (5) to the form
\begin{equation}
(Z^t)^{-1}q=V^{-1}ZV,\,\,\,{\rm tr}\,Z=1+q.
\end{equation}

The family of all solutions to the QYBE over the field $K= {\bf C}$ is
described by the following
\begin{proposition}
{\em \cite{G}} The pair $(Z,\,V)$ is a solution of the system {\em (6)} iff
the matrix $Z$ is such that ${\rm tr}\,Z=1+q$ and its Jordan form contains
along with any cell corresponding to an eigenvalue $x$, an analogous cell
with eigenvalue $q/x$ {\em (}with the same multiplicity{\em )}.
\end{proposition}
\begin{remark}
{\em Let us note that $U$ and $V$ are transformed under changes of base as
bilinear form matrices, while $Z$ is transformed as an operator matrix (their
transformations are coordinated and the relations (6) are stable). So,
assuming $K= {\bf C}$, we can represent the operator $Z$ in Jordan form by an
appropriate choice of base.  Moreover, we can assume that the cells with
eigenvalues $x$ and $q/x$ are in positions symmetric to each other with
respect to the center of the matrix $Z$. Observe that if the number of the
cells is odd the eigenvalue of the middle one is $\pm\sqrt q$.}
\end{remark}

It is not difficult to see that for such a choice of base the tensor $V$ can
be taken in the form of a skew-diagonal matrix (i.e., possessing nontrivial
terms only at the auxiliary diagonal). Let us fix such a matrix $V_0$ and
note that all other $V$ satisfying (6) are of the form $V=WV_0$, where $W$
commutes with $Z$. In the sequel we assume that a base possessing these
properties is fixed.

Let us observe that in case under consideration we have $M=UV,\,\, N=VU$.
Moreover, relations (4) take the form $m^2=q(1+q)^{-2}$. Using the relation
$U^t=(1+q)^{-1}\,V^{-1}\,Z$, we can transform the equality $UV=m\,{\rm id}$
to
\begin{equation}
Z=(1+q)m V(V^t)^{-1}.
\end{equation}

Let us assume that $Z$ has a simple spectrum, i.e., its eigenvalues are
pairwise distinct. So, its Jordan form is diagonal: $Z={\rm diag}(z_1,\dots,
z_n)$. The family of diagonal $Z$, satisfying conditions of Proposition 2 and
fulfilling the only relation ${\rm tr}\,Z=1+q$, can be parametrized by $(z_1,
\dots ,z_r)$ with $r=n/2$, if $n$ is even, and with $r=(n-1)/2$, if $n$ is
odd (if $n$ is odd we have also a choice for the value of $z_{r+1}=\pm\sqrt
q$).

Since $Z$ has a simple spectrum, any $W$ commuting with $Z$ is also diagonal
(with arbitrary diagonal entries). This implies that $V$ satisfies (6) iff it
is skew-diagonal with arbitrary entries at the auxiliary diagonal. Therefore
$U$ is also skew-symmetric. Thus, we have $v^{ij}\not=0,\, u_{ij}\not=0$ iff
$i+j=n$. For the sake of simplicity, we will use the notation $v^i\,\,(u_i)$
instead of $v^{i\,\,n+1-i}\,\,(u_{i\,\,n+1-i})$. Let us note that $z_i=(1+q)
u_iv^i$ (up to the end of this section there is no summation over repeated
indices).

It is easy to see that relation (7) is satisfied iff the entries $v_i$
fulfill the system
\begin{equation}
m(1+q)v^i/v^{n-i+1}=z_i,\,\,1\leq i\leq n.
\end{equation}
This system is consistent by virtue of the relations
$$z_i\,z_{n-i+1}=q,\,\,m^2=q(1+q)^{-2}.$$
Moreover, the family of the solutions of the system (8) can be parametrized
by $(v_1,\dots ,v_{r})$. Let us note that if $n$ is odd, the value of $z_{r
+1}=\pm{\sqrt q}$ depends on that of $m=\pm{\sqrt q}(1+q)^{-1}$, namely, we
have $z_{r+1}=m(1+q)$.

Thus, we have proved the following
\begin{proposition} Let $K= {\bf C}$, $S$ be a Hecke symmetry of TL type and
$Z$ be the corresponding tensor described in Proposition 2 with a simple
spectrum {\em (}a parametrization of all such tensors $Z$ was given above{\em
)}. Then the dual quasitriangular structure defined on the algebra $\ar$
can be descended on the quantum cogroup $\Fun(SL(S))$ iff $V$ is a
skew-diagonal matrix with the entries $v_{i\,n+1-i}=v_i$ satisfying the
system {\em (8)}. This system is always compatible and the family of its
solutions can be parametrized as above.
\end{proposition}
\begin{remark}{\em
Let us observe that the Hecke symmetries of TL type such that the operator $U
V$ is scalar are just those introduced in \cite{DL} (the authors of \cite{DL}
use another normalization of the operator $S$).}
\end{remark}

\section{Nondegeneracy of the canonical pairing}
In the present section we show that when $n={\rm dim}\,{\bf V}>2$ and $q$ is
generic, the canonical pairing $\llangle\,\,,\,\,\rrangle _c$ is
nondegenerate for those even Hecke symmetries of TL type whose operator $Z$
has a generic simple spectrum. As above, we assume that $Z$ has a diagonal
form in a chosen base and therefore the tensors $U$ and $V$ are
skew-diagonal. In the sequel we put $c=1$.

Thus we have the pairing $\llangle t_i^j,\, t_k^l\rrangle=S_{ki}^{jl}$. To
show that it is nondegenerate we will compute the Gram determinant, i.e., the
determinant of the Gram matrix. The rows and the columns of this matrix are
labeled by the bi-index $(i,j)$ running over the set
$$(1,1),\dots ,(1,n),(2,1),\dots ,(2,n),\dots ,(n,1),\dots ,(n,n).$$
So, the term $\llangle t_i^j,t_k^l\rrangle=R_{ik}^{jl}=S_{ki}^{jl}$ is
situated at the intersection of the $(i,j)$-row and the $(k,l)$-column.

Let us note that if $S$ is a Hecke symmetry of TL type, then all the
entries of the matrix $S_{ki}^{jl}$ are equal to zero unless either $i=j,\,
k=l$ or $i+k=j+l=n+1$. So, we have just two nonzero elements in the
$(i,j)$-row namely, $R_{ij}^{ji}$ and $R_{i\,n+1-i}^{j\,n+1-j}$, if $i+
j\not=n+1$, and only one, namely, $R_{i\,n+1-i}^{n+1-i\,i}$, if $i+j=n+1$. A
similar statement is valid for the columns.

This yields that if $i+j\not=n+1$ then the $(i,j)$- and $(n+1-j,n+1-i)$-rows
and $(j,i)$- and the $(n+1-i,n+1-j)$-columns possess just four nontrivial
elements
$$R_{ij}^{ji},\, R_{i\,n+1-i}^{j\,n+1-j},\,
R_{n+1-j\,j}^{n+1-i\,i},\, R_{n+1-j\,n+1-i}^{n+1-i\,n+1-j}$$
situated at their intersections. If $i+j=n+1$, then two rows (columns) are
merged into one and the only nontrivial element $R_{i\,n+1-i}^{n+1-i\,i}$
belongs to the intersection of the $(i,\, n+1-i)$-row and the $(n+1-i,
\,i)$-column.

For example, for $n=3$ we have the following Gram matrix
$${\bf G}=\left(\begin{array}{ccccccccc}
R_{11}^{11} & 0 & 0 & 0 & 0 & 0 & 0 & 0 & R_{13}^{13}\\
0 & 0 & 0 & R^{21}_{12} & 0 & 0 & 0 & R^{22}_{13} & 0\\
0 & 0 & 0 & 0 & 0 & 0 & R^{31}_{13} & 0 & 0\\
0 & R^{12}_{21} & 0 & 0 & 0 & R^{13}_{22} & 0 & 0 & 0\\
0 & 0 & 0 & 0 & R_{22}^{22} & 0 & 0 & 0 & 0\\
0 & 0 & 0 & R^{31}_{22} & 0 & 0 & 0 & R^{32}_{23} & 0\\
0 & 0 & R^{13}_{31} & 0 & 0 & 0 & 0 & 0 & 0\\
0 & R^{22}_{31} & 0 & 0 & 0 & R^{23}_{32} & 0 & 0 & 0\\
R_{31}^{31} & 0 & 0 & 0 & 0 & 0 & 0 & 0 & R_{33}^{33}
\end{array}.\right)$$

By changing the order of the rows and columns, we can reduce this matrix
to a block-diagonal form, where all blocks are either one-dimensional and
consist of elements $R_{i\,n+1-i}^{n+1-i\,i}$ or two-dimensional and have the
following form
$$\left(\begin{array}{cc}
R_{ij}^{ji} & R_{i\,n+1-i}^{j\,n+1-j}\\
R_{n+1-j\,j}^{n+1-i\,i} & R_{n+1-j\,n+1-i}^{n+1-i\,n+1-j}
\end{array}\right).$$
Denote by $I(n)$ the set of indices $i,j$ satisfying $1\leq i,j\leq n,i+j\neq
n+1$. Then $\det\,{\bf G}$ is equal up to a sign to
$$\prod_{1\leq i\leq n}R_{i\,n+1-i}^{n+1-i\,i}\sqrt{\left|\prod_{I(n)}
\det\left(\begin{array}{cc}
R_{ij}^{ji} & R_{i\,n+1-i}^{j\,n+1-j}\\
R_{n+1-j\,j}^{n+1-i\,i} & R_{n+1-j\,n+1-i}^{n+1-i\,n+1-j}
\end{array}\right)\right|}.$$
(The root is motivated by the fact that any factor in the second product is
taken two times.)

By straightforward calculations it is not difficult to see that
$$\prod_{1\leq i\leq n}R_{i\,n+1-i}^{n+1-i\,i}=\prod_{1\leq i\leq n}(q-z_i)$$
and
$$\prod_{I(n)}\det\left(\begin{array}{cc}
R_{ij}^{ji} & R_{i\,n+1-i}^{j\,n+1-j}\\
R_{n+1-j\,j}^{n+1-i\,i} & R_{n+1-j\,n+1-i}^{n+1-i\,n+1-j}
\end{array}\right)=\prod_{I(n)}(q^2-z_{n+1-i}z_j).$$
Thus, we have proven the following
\begin{proposition}
$$(\det\,{\bf G})^2=\prod_{1\leq i\leq n} (q-z_i)^2
\prod_{I(n)}(q^2-z_{n+1-i}z_j).$$
\end{proposition}

It is interesting to observe that the final expression depends only on the
matrix $Z$. This enables us to state the following
\begin{proposition} Let us assume that $S:{\bf V}^{\ot 2}\to {\bf V}^{\ot
2}$ is an even Hecke symmetry of TL type, the corresponding operator $Z$
possesses a simple spectrum and $n={\rm dim}\, {\bf V}\geq 4$. Then for a
generic $q$ and for a generic $Z$ {\em (}of such type{\em )} the canonical
pairing is nondegenerate.
\end{proposition}

{\it Proof.} The set where the determinant $\det\, {\bf G}$ vanishes is an
algebraic variety in the space ${\bf C}^{n+1}$ generated by the
indeterminates $z_i$ and $q$. It suffices to show that this variety is not
contained in that defined by
$$z_i\,z_{n+1-i}=q,\ \ \sum z_i=1+q.$$
Let us decompose the expression
$$\prod_{I(n)} (q^2-z_{n+1-i}z_j) $$
into the product of factors with $i=j$ and those with $i\not=j$.

If $i=j$, we have $z_{n+1-i}z_j=q$. Thus, the above product is equal to
$$(q^2-q)^n\prod_{I(n),\,i\not= j}(q^2-z_{n+1-i}z_j).$$
So, for $q$ such that $q\not=0,q\not=1$, we have $\det\,{\bf G}=0$ iff
$z_i=q$ for some $i$ or $q^2=z_{n+1-i}z_j$ for some $i\not=j$.

From the above parametrization it is evident that if $n\geq 4$ there exists a
matrix $Z={\rm diag}(z_1,\dots ,z_n)$ satisfying the conditions of
Proposition 2 and such that $\det\,{\bf G}\not=0$.

Let us consider the case $n=3$ separately. In this case the set of such
diagonal matrices $Z$ is parametrized by $z_1$ (after choosing a value of
$z_2=\pm\sqrt q$) satisfying the equation $z_1\pm\sqrt q+q/z_1=1+q$. This
equation has two solutions for any choice of the value $\pm\sqrt q$. It is
not difficult to see that for a generic $q$ we have $\det\,{\bf G}\not=0$ for
any of these four values of $z_1$. \ $\Box$

Let us note that the canonical pairing becomes degenerate if $q=1$ (this case
corresponds to an involutive symmetry $(S^2={\rm id})$) or if $n=2$, since in
this case the system $z_1z_2=q,\, z_1+z_2=1+q$ has two solutions $z_1=1,\,
z_2=q$ and $z_1=q,\,z_2=1$ for which the product $(q-z_1)(q-z_2)$ vanishes.
And we always have $\det\,{\bf G}=0$.

This is the principal reason why the latter case $(n=2)$ which corresponds to
a quasiclassical Yang--Baxter operator $S$ (it is in fact the only
quasiclassical case related to the TL algebra) differs crucially from the
nonquasiclassical ones $(n>2)$.

Thus, according to the above construction, we can convert the right
$\Fun(SL(S))$-comodules ${\bf V}^{\ot m}$ into left
$\Fun(SL(S))^{op}$-modules (assuming $S$ to be an even Hecke symmetry of TL
type, $q$ to be generic and $Z$ to have a simple spectrum, also generic) and
therefore into right $\Fun(SL(S))$-modules.

\section{Discussion of a  possible representation theory}
Our further aim is to construct some representation theory of the
algebra $\Fun(SL(S))$ equipped with the above action. Conjecturally, it looks
like that of $SL(2)$. Let us denote $V_m$ the symmetric component of ${\bf V}
^{\ot m}$. (Let us note that in classical and quasiclassical cases $m/2$ is
just the spin of the representation $U(sl(2))\to{\rm End}(V_m)$ or $U_q
(sl(2))\to{\rm End}(V_m)$.)

It seems very plausible that similarly to the $SL(2)$- or $U_q(sl(2))$-case
we have for a generic $q$ the following properties
\begin{itemize}
\item the $\Fun(SL(S))$-modules $V_m$ are irreducible,
\item any irreducible finite-dimensional $\Fun(SL(S))$-module is isomorphic
to one of $V_m$,
\item any finite-dimensional $\Fun(SL(S))$-module is completely reducible,
\item we have the classical formula $V_i\ot V_j=\oplus_{\vert i-j\vert\leq
k\leq i+j}V_k$.
\end{itemize}
To motivate the latter formula, let us show that it is satisfied at least
``in sense of dimensions", i.e.,
\begin{equation}
{\rm dim}\,V_i\ot {\rm dim}\,V_j=\sum_{\vert i-j\vert\leq k\leq i+j}
{\rm dim}\,V_k.
\end{equation}
Indeed, using the fact that the Poincar\'e series of the symmetric algebra of
the space ${\bf V}$ is equal to $(t^2-nt+1)^{-1}$, one can see that
$${\rm dim}\,V_i=\alpha^i+\alpha^{i-2}+\alpha^{i-4}+\dots+\alpha^{-i},$$
where $\alpha=n/2+\sqrt{(n/2)^2-1}$ is a root of the equation $t^2-nt+1=0$.
Then relation (9) can be established by straightforward calculations. The
details are left to the reader.

\end{document}